\documentclass[12pt]{article}

\usepackage{amsmath,empheq}
\usepackage{geometry}
\usepackage{graphicx}
\usepackage{subcaption}
\usepackage{placeins}
\usepackage{seqsplit}
\usepackage{physics}
\usepackage{mhchem}
\usepackage{subcaption}
\usepackage{mathptmx}
\usepackage{setspace}
\usepackage{enumitem}
\usepackage{comment}
\usepackage{titlesec}
\usepackage{float}
\usepackage{mathrsfs}
\usepackage{listings}
\usepackage{siunitx}
\usepackage{color}
\usepackage{pythonhighlight}
\usepackage[symbol]{footmisc}

\definecolor{dkgreen}{rgb}{0,0.6,0}
\definecolor{gray}{rgb}{0.5,0.5,0.5}
\definecolor{mauve}{rgb}{0.58,0,0.82}

\usepackage{todonotes}
\usepackage{amssymb}
\usepackage[numbers]{natbib}
\usepackage{tikz}
\usetikzlibrary{patterns}
\usetikzlibrary{shapes.geometric, arrows}
\usetikzlibrary{decorations.pathreplacing}

\usepackage{pgfplots}

\usepackage[dotinlabels]{titletoc}
\numberwithin{equation}{section}
\titleformat{\section}[hang]{\bfseries\large}{\thesection.}{0.4em}{}
\titlelabel{\thetitle.\quad}
\graphicspath{{.}{figures/}}

\title{\textbf{PeTTO: Leveraging GPUs to Accelerate Topology Optimization with the Pseudo-Transient Methods}}
\author{Mingyuan Yang$^{1}$ \and Qian Yu$^{1}$ \and Chao Yang$^{1, *}$}
\date{$^1$ School of Mathematical Sciences, Peking University, Beijing 100871, PR China}

\begin{document}
\maketitle

\begin{sloppypar}

\section*{Abstract}

We present a Pseudo-Transient Topology Optimization (PeTTO) approach that can leverage graphics processing units (GPUs) to efficiently solve single-material and multi-material topology optimization problems. 
By integrating PeTTO with phase field methods, the partial differential equations (PDEs) constrained optimization problem in topology optimization is transformed into a set of time dependent PDEs, which can be analyzed using the knowledge of transient physics.
The sensitivities with respect to the design variable are calculated with the automatic differentiation which help avoid tedious and error-prone manual derivations.  
The overall system of equations is efficiently solved using a hybrid of the pseudo-transient method and the accelerated pseudo-transient method, balancing the convergence rate and numerical stability.
A variety of numerical examples are presented to demonstrate the effectiveness and efficiency of the proposed PeTTO approach. 
These examples cover different physics scenarios including mechanical and thermal problems, as well as single-material and multi-materials cases in both 2D and 3D.
The numerical results show a 40- to 50-fold speedup when running the same PeTTO code on a single GPU compared to desktop CPUs.
This work helps bridge the gap between high-performance computing and topology optimization, potentially enabling faster and better designs for real-world problems. 

\vspace{1em}
  
\noindent Keywords: Topology optimization; Pseudo-transient method; Phase field method; Differential programming; High performance computing

\footnotetext{*Corresponding author.}
\footnotetext{ Email addresses: mingyuanyang@pku.edu.cn (Mingyuan Yang), qianyu@pku.edu.cn (Qian Yu), chao\_yang@pku.edu.cn (Chao Yang).}

\newpage

\section{Introduction}\label{introduction}

Topology optimization is a powerful numerical technique for determining the optimal distribution of material within a given space, subject to specific design constraints and performance requirements.
Since its emergence in 1989, topology optimization has been extensively studied in academia and widely adopted in various industries. 
In the paper by Yap et al.~\cite{yap2023topology}, topology optimization is employed to enhance the design of micro-unmanned aerial vehicles through 3D printing, resulting in significant weight reduction and increased reliability.  
In the paper by Kendibilir et al.~\cite{kendibilir2022peridynamics}, the authors combine topology optimization together with peridynamics to optimize structures in the additive manufacturing process considering surface cracks. 
In the paper by Zhang et al.~\cite{zhang2021paved}, topology optimization is used to design paved pedestrian guideway networks, balancing traffic congestion capability and infrastructure construction costs.
The growing adopting of topology optimization across various industries is driving up the demand for large-scale simulations~\cite{mukherjee2021accelerating,herrero2023parallel}.

The acceleration of topology optimization is primarily achieved through three strategies, each one targeting the reduction of the computational time required for a specific step in the topology optimization process.
The three strategies are: solving the governing partial differential equations (PDEs) faster, improving the efficiency of sensitivity calculations, and updating the design variables in a more computationally effective manner.  
Solving the governing PDEs is computationally intensive and accounts for a significant portion of the total time required for topology optimization.
Therefore, improving the efficiency of this step is crucial for accelerating the overall optimization process~\cite{mukherjee2021accelerating,aage2015topology}. 
In the paper by Perez et al.~\cite{herrero2023parallel}, the authors present an efficient parallel geometric multigrid~(GMG) implementation for preconditioning Krylov methods to solve the governing PDEs.
In addition, the authors use adaptive meshing techniques to reduce the number of computational elements where fine structures are not needed.
The method is implemented using the PETSc library for parallel computing on multiple CPU cores, significantly accelerating the solving process~\cite{herrero2023parallel}.
Besides multi-CPU computing schemes, graphics processing units (GPUs) have also advanced in scientific computing fields in recent years due to their extraordinary parallel computing capabilities~\cite{herrero2021multi,kazakis2017topology,challis2014high}.
In the paper by Perez et al.~\cite{herrero2021multi}, the authors implement an algebraic multigrid~(AMG) preconditioned conjugate gradient (PCG) solver on GPUs to solve the linear system of equations arising from topology optimization.
The authors highlight the advantages of using GPUs for their algorithms, noting reduced memory requirements and improved computing efficiency.

Although significant computational speedup has been achieved using GPUs, traditional Krylov methods may not fully exploit the potentials of GPUs~\cite{rass2022assessing,wang2022physics}.
In their paper, Rass et al.~\cite{rass2022assessing} propose the accelerated pseudo-transient~(PT) method, which enhances the first-order PT method by incorporating higher-order derivatives with respect to pseudo time.
Compared to Krylov-type methods for solving linear system of equations, this method updates each grid point locally, thereby eliminating the need for global communication~\cite{rass2022assessing}.
In the paper by Wang et al.~\cite{wang2022physics}, the authors employ the accelerated PT method to model geological processes coupling solid deformation and fluid migration. 
These progress have inspired us to explore the effectiveness of the accelerated PT method in solving the PDEs that arise in topology optimization problems. 
Details and considerations of our implementation will be discussed in Section~\ref{model_formulation}.

Sensitivity calculation is an essential step in topology optimization to guide the update of design variables~\cite{yu2021first,yu2023second,li2023convolution}.
Manual derivation of sensitivities has been a traditional approach in topology optimization, however, it is time-consuming, error-prone, and becomes increasingly complex as more constraints are added to meet realistic requirements~\cite{chandrasekhar2021auto,he2022multiphysics,mowlavi2023topology}.
Automatic differentiation~(AD) has gained attention in the topology optimization field as an efficient alternative to manual derivation for calculating sensitivities~\cite{mallon2023neural,chandrasekhar2021auto,he2022multiphysics}. 
In the paper by Chandrasekhar et al.~\cite{chandrasekhar2021auto}, the authors demonstrate that the computational cost of calculating sensitivities is nearly equivalent to that of using analytical formulas, with the added advantage that only the forward evaluation needs to be defined. 
In the paper by Mallon et al.~\cite{mallon2023neural}, it demonstrates the seamless integration of AD step with other computational steps, enabling the possibility of efficient inverse modeling.
In this work, we implement our algorithm using JAX~\cite{jax2018github}, a Python library for high-performance numerical computing and machine learning~\cite{bezgin2023jax,xue2023jax,mistani2023jax}. JAX enables automatic differentiation through Python functions which ease the algorithm implementation. 
Details of sensitivity calculations can be found in~\cite{chandrasekhar2021auto}, 
our implementation is similar to the code shown in~\cite{chandrasekhar2021auto} and will be available upon request. 

Numerous topology optimization methods have been developed for various fields, including the homogenization-based approach~\cite{xia2015design}, the solid isotropic material with penalization method~(SIMP)~\cite{da2022some,liu2014efficient}, 
the level set-based methods~(LS)~\cite{mallon2023neural} and the phase-field based methods~(PF)~\cite{zhou2007multimaterial,yu2021first,yu2023second}, among others.
The phase field method, originally developed to describe phase transition phenomena, has been successful applied in the field of topology optimization~\cite{dede2012isogeometric,carraturo2021additive}.
The phase field method permits hole nucleation during the optimization process, resulting in similar final structures as other methods.
Moreover, its diffusive interface representation eliminates the need for sequential variable adjustments and inherently prevents the checkerboard effects~\cite{dede2012isogeometric}.
In this work, we utilize the Cahn-Hilliard~(CH) equation to evolve the topological changes in the structures. 
The explicit scheme of the CH equation, also implemented on GPUs, integrates seamlessly with the PT method to effectively and efficiently obtain the final optimized structure. 
To the best of our knowledge, this is the first time that a topology optimization problem has been reformulated as a set of time-dependent PDEs coupled through sensitivities.
Our approach, PeTTO, implemented on GPUs and leveraging automatic differentiation for efficient sensitivity calculation, complements existing methods for topology optimization and has potential to lead to superior designs in various fields. 

The remainder of this paper is organized as follows. Section~\ref{model_formulation} provides a brief introduction to the methods used in this work and presents the coupled formulation for the topology optimization problems considered herein.
Section~\ref{numerical_examples} presents a range of numerical examples that demonstrate the effectiveness and the efficiency of the proposed approach, from 2D benchmark problems to 3D realistic applications.   
Section~\ref{conclusion} concludes this work and outlines several promising directions for the future research.

\section{Model Formulation}\label{model_formulation}

\subsection{PDE constrained topology optimization}

A general PDE constrained topology optimization problem can be stated as follow:
\begin{align}
    \min\limits_{\phi} ~~ \mathscr{J} \left(u \left( \phi\right), \phi\right), \nonumber
    \\
    s.t. ~ \mathscr{R} \left(u\left( \phi\right), \phi\right) = 0,
    \\
    \mathscr{C}\left( \phi\right) = 0, \nonumber
\end{align}
where $\phi$ is the design variable of interest, $u$ is the solution to the governing stationary PDE, $\mathscr{J}$ is the objective function,
$\mathscr{R}$ is the PDE residual and $\mathscr{C}$ is the constraint function of $\phi$. 
The governing PDE describes the underlying physical behavior of the system. The goal is to find the optimal $\phi$ that minimizes the objective function subject to the constraints imposed by the PDE. 

Figure~\ref{a typical workflow of topology optimization} shows a typical workflow for PDE constrained topology optimization.  
The process starts with an initial distribution of the design variable $\phi$ across the computational domain.
The governing PDE, such as the elasticity equation or the Poisson equation, is solved with respect to this design using a suitable numerical method, such as the finite element method.
An optimization method is then used to minimize the objective function which is defined in terms of the design variable $\phi$ and the solution $u$ to the PDE. 
The convergence criteria, such as the change in the objective function or the design variable, are checked against to determine if the optimization continues.
If all the stopping criteria are met, the final design is returned, Otherwise, the variable $\phi$ is updated and the process is repeated for another optimization loop. 

In the following sections, we will provide a concise overview of the methods our approach employs to tackle a topology optimization problem.

\tikzstyle{terminator} = [rectangle, draw, text centered, rounded corners, minimum height=2em]
\tikzstyle{process} = [rectangle, draw, text centered, minimum height=2em]
\tikzstyle{decision} = [diamond, draw, text centered, minimum height=2em]
\tikzstyle{data}=[trapezium, draw, text centered, trapezium left angle=60, trapezium right angle=120, minimum height=2em]
\tikzstyle{connector} = [draw, -latex']

\begin{figure}[hbt!]
\centering
    \begin{tikzpicture}
        \node [terminator, fill=blue!20] at (0,0) (start) {\textbf{Start}};
        \node [data, fill=blue!20] at (0,-2) (theta) {$\phi$};
        \node [data, fill=blue!20] at (0,-4) (constraint) {$\mathscr{C}\left( \phi\right) = 0$};
        \node [data, fill=blue!20] at (0,-6) (PDE) {$\mathscr{R} \left(u\left(\phi\right), \phi\right) = 0$};
        \node [data, fill=blue!20] at (0,-8) (u) {$u$};
        \node [process, fill=red!20] at (3.5,-5) (update) {Update $\phi$};
        \node [process, fill=red!20] at (0,-10) (success) {$\min\limits_{\phi} \mathscr{J} ?$};
        \node [terminator, fill=blue!20] at (0,-12) (end) {\textbf{End}};

        \node[draw=none] at (1.85, -9.7) (no) {No};
        \node[draw=none] at (0.35, -11) (yes) {Yes};
        \path [connector] (start) -- (theta);
        \path [connector] (theta) -- (constraint);
        \path [connector] (constraint) -- (PDE);
        \path [connector] (PDE) -- (u);
        \path [connector] (u) -- (success);
        \path [connector] (success) -| (update);
        \path [connector] (update) |- (theta);
        \path [connector] (success) -- (end);
    \end{tikzpicture}
    \caption{A typical schematic of workflow for PDE constrained topology optimization.}\label{a typical workflow of topology optimization}
\end{figure}
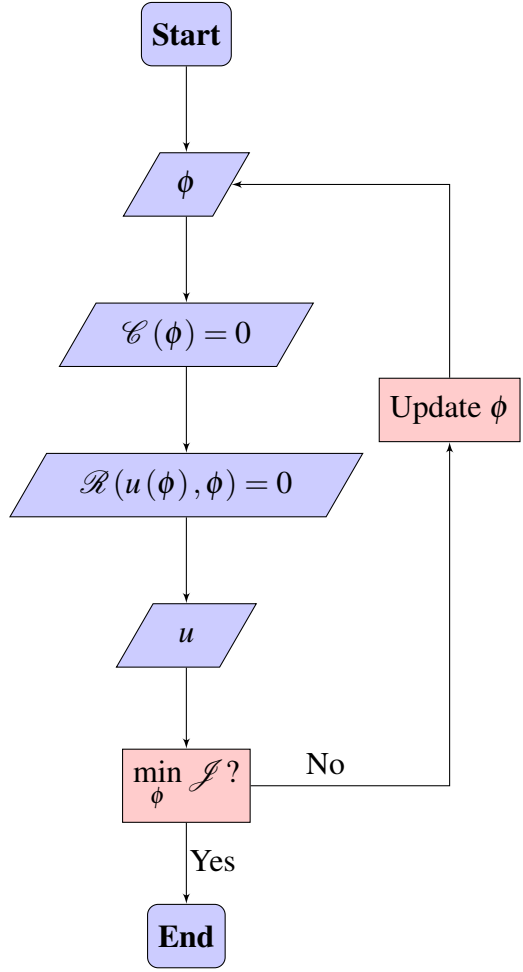

\FloatBarrier 

\subsection{The pseudo transient method}\label{the pseudo transient method}

The pseudo transient method, also known as the pseudo transient continuation, is a physics-motivated iterative method to solve a wide range of PDEs.
The core idea of PT is to transform a stationary PDE problem into a transient PDE problem where standard numerical methods for time-dependent PDEs apply. 
It has a long history in the literature and can be dated back to 1950s~\cite{frankel1950convergence}.
Consider the Poisson equation in the context of a heat conduction topology optimization problem:

\begin{equation}
    - \nabla \cdot \left( \kappa \nabla T \right) = f ~~~ in ~ \Omega,\label{Poisson equation}
\end{equation}
$\Omega$ is the computational domain, $\kappa$ is the thermal conductivity tensor of the material, $T$ is the temperature field of the domain and $f$ is the thermal source.
If we add a time-dependent term into Equation~\eqref{Poisson equation}, it becomes a heat equation, 

\begin{equation}
    \frac{\partial T}{\partial t} = \nabla \cdot \left(\kappa \nabla T\right) + f,\label{heat equation}
\end{equation}
$t$ is referred as pseudo-time here as we are not interested in the solution of $T$ at a particular time of $t$. 
The solution to Equation~\eqref{heat equation} can be achieved as $t$ goes to infinity.
In practice, we perform the iteration in an explicit fashion to fully exploit the power of GPUs as discussed in~\cite{rass2022assessing} and the discrete formulation writes as:

\begin{equation}
    \frac{T^{n+1} - T^{n}}{\Delta t} = \nabla \cdot \left(\kappa \nabla T^{n}\right) + f,\label{discretized heat equation}
\end{equation}
where $T^{n}$ and $T^{n+1}$ are the solutions at current and next time step, $\Delta t$ is the pseudo time step. Though it converges to the steady state solution if the selection of $\Delta t$ satisfies the CFL condition, this formulation suffers from slow convergence rate. 
If we consider having $n_x$ grids in each spatial dimension, the number of steps in pseudo-time to reach convergence is proportional to $n_x^{2}$~\cite{rass2022assessing,wang2022physics}. 
One remedy is to use the accelerated pseudo transient method proposed in~\cite{rass2022assessing}, further including a second order derivative term with respect to the pseudo-time $t$. 
The Equation~\eqref{Poisson equation} can be reformulated as:

\begin{equation}
    \theta \frac{\partial^2 T}{\partial t^2} + \frac{\partial T}{\partial t} = \nabla \cdot \left(\kappa \nabla T\right) + f,\label{damped wave equation}
\end{equation}
which is a damped wave equation and $\theta$ is the damping factor. 
Again, the explicit discrete formulation can be written as:
\begin{equation}
    \theta \frac{T^{n+1} - 2T^{n} + T^{n-1}}{\Delta t^2} + \frac{T^{n} - T^{n-1}}{\Delta t} = \nabla \cdot \left(\kappa \nabla T^{n}\right) + f,\label{discretized damped wave equation}
\end{equation}
where $T^{n-1}$, $T^{n}$ and $T^{n+1}$ are solutions at three consecutive time steps.
The selections of $\theta$ and $\Delta t$ significantly influence the performance and the stability of the accelerated PT method. For this Poisson equation case, the iterations required to reach convergence can be reduced to being proportional to $n_x$. 
Our numerical experiments showed that, while the accelerated PT method can effectively update the solution to the governing PDEs, it can occasionally lead to numerical instability.
To balance computational efficiency and stability, we propose a hybrid approach that combines the PT method and the accelerated PT method within a single optimization loop.
In our numerical examples, detailed in Section~\ref{numerical_examples}, we employ a fixed number of accelerated PT iterations with a subsequent fixed number of PT iterations to update the solution to governing PDEs.

\subsection{Sensitivity calculation}\label{sensitivity calculation}

A critical aspect of topology optimization is updating the design variable using sensitivities~\cite{chandrasekhar2021auto}. 
To provide a clearer understanding, let us once again use the heat conduction problem as an illustrative example: 
\begin{gather}
    J_h = \int_{\Omega} \nabla T^T \kappa \nabla T d\Omega,\label{thermal compliance}
    \\
    \kappa = \sum_i \kappa_i \phi_i,~~i=1,2,...,N, 
\end{gather}
where $J_h$ is the thermal compliance, calculated using the gradient of temperature $\nabla T$ and the overall thermal conductivity tensor $\kappa$.
$\kappa$ is further calculated using the volumetric fraction of different materials $\phi_i$ and their respective conductivity $\kappa_i$.
$J_h$ is a measure of the system's ability to dissipate heat. In a typical heat conductivity topology optimization problem, the goal is to minimize this quantity. 
To achieve this, we need to calculate the partial derivatives of $J_h$ with respect to the design variable $\phi_i$, and then update $\phi_i$:  
\begin{equation}
    \phi_i^{n+1} = \phi_i^{n} - \alpha \frac{\partial J_h}{\partial \phi_i^{n}},~~i=1,2,...,N.
\end{equation}
Here $i=1,2,...,N$ and $N$ is the number of materials. 
\FloatBarrier 

\begin{figure}[hbt!]
    \begin{python}
        #heat compliance definition 
        #k -- thermal conductivity tensor
        #T -- temperature field
        def heat_compliance(k, T, dx, dy):
            dTdx, dTdy = jnp.gradient(T, dx, dy)
            J = jnp.sum(k *  (dTdx**2 +  dTdy**2))
            return J
    
        #sensitivity 
        dJdk = jax.grad(heat_compliance)(k, T, dx, dy)
    
        #update
        k -= alpha *  dJdk
    \end{python}
\caption{JAX code for heat compliance definition, sensitivity calculation using AD and the thermal conductivity update.}\label{heat sensitivity code}
\end{figure}


In this work, we leverage the automatic differentiation technique to calculate sensitivities. As shown in Figure~\ref{heat sensitivity code},
only minimal additional effort is required to compute the sensitivities using the forward evaluation function while maintaining efficiency. 
Moreover, this approach can be extended to arbitrarily complex sensitivity calculation to meet various design specifications, including but not limited to compliance minimization,
stress minimization, volume control, size control~\cite{chandrasekhar2021auto}.

Another key constraint in topology optimization is the volume constraint, which ensures that the final structure remains lightweight without significantly compromising performance. 
The objective function with respect to the volume constraint can be defined as:
\begin{equation}
    J_v = \sum_i \left(\int_{\Omega} \phi_i d\Omega - V_i\right)^2,
\end{equation}
$\phi_i$ is the density of each phase across the domain and $V_i$ is the respective target volume fraction.
Additionally, in multi-material topology optimization, there is a partition of unity constraint requiring that the densities of all phases must sum up to one at every point in the domain:
\begin{equation}
    J_1 = \int_{\Omega}\left(\sum_i \phi_i - 1\right)^2 d\Omega.
\end{equation}

Therefore, the overall objective function for sensitivity calculation of the heat conduction topology optimization problem can be written as:
\begin{equation}
    J = \alpha_h J_h + \alpha_v J_v + \alpha_1 J_1,
\end{equation}
where $\alpha_h$, $\alpha_v$ and $\alpha_1$ are weights to be selected to balance each objective quantity. 
The design variable can be updated as:
\begin{equation}
    \phi_i^{n+1} = \phi_i^{n} - \left(\alpha_h \frac{\partial J_h}{\partial \phi_i^{n}} + \alpha_v \frac{\partial J_v}{\partial \phi_i^{n}} + \alpha_1 \frac{\partial J_1}{\partial \phi_i^{n}}\right).
\end{equation}

\subsection{Topology optimization with the phase field method}
In this work, we utilize the Cahn-Hilliard model to drive the topological changes in structures.
A typical Cahn-Hilliard model can write as:
\begin{gather}
    \begin{cases}\label{cahn hilliard model}
        \frac{\partial \phi}{\partial t} = D \nabla^2 \mu,~~in~\Omega, 
        \\
        \\
        \mu = \frac{\partial w\left(\phi\right)}{\partial \phi} - \gamma \nabla^2 \phi,
        \\
        \\
        \nabla \mu = 0,~~in~d\Omega, 
        \\
        \\
        \nabla \phi = 0,~~in~d\Omega,
    \end{cases}
\end{gather}
where $\phi$ is an order parameter in the phase field theory, D denotes the diffusion coefficient, $w\left(\phi\right)$ represents the double well potential function and $\gamma$ is the interfacial thickness parameter.  
Here we define $\phi$ to range between 0 and 1, representing the density as it appears in the sensitivity calculation, with $w\left(\phi\right)$ having local minima at 0 and 1. By evolving the scalar variable $\phi$ using the Cahn-Hilliard model, it promotes phase separation while preserving smooth interfacial transitions~\cite{zhou2007multimaterial}. 
The no-flux boundary condition for both $\mu$ and $\phi$ further ensures mass conservation throughout the phase evolution.
Here, $w\left(\phi\right)$ is defined as:
\begin{equation}
    w(\phi) = \frac{1}{64} \sin^2(\phi \pi).
\end{equation}

\FloatBarrier 

\begin{figure}[hbt!]
    \centering
    \begin{subfigure}{0.45\textwidth}
        \centering
        \includegraphics[width=0.9\linewidth]{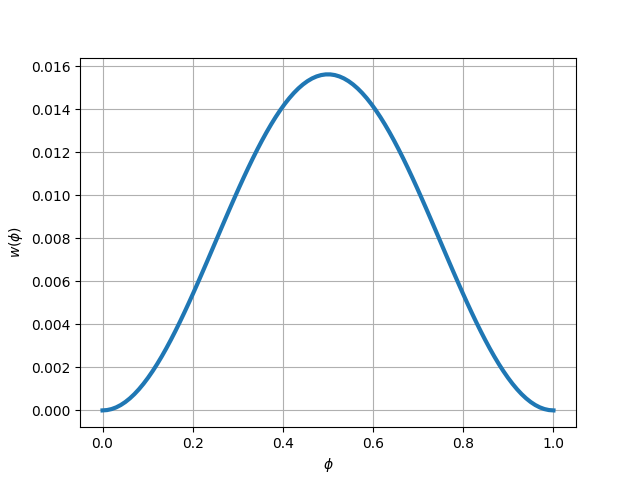}
        \caption{$w\left(\phi\right)$ vs. $\phi$.}
    \end{subfigure}
    \begin{subfigure}{0.45\textwidth}
        \centering
        \includegraphics[width=0.9\linewidth]{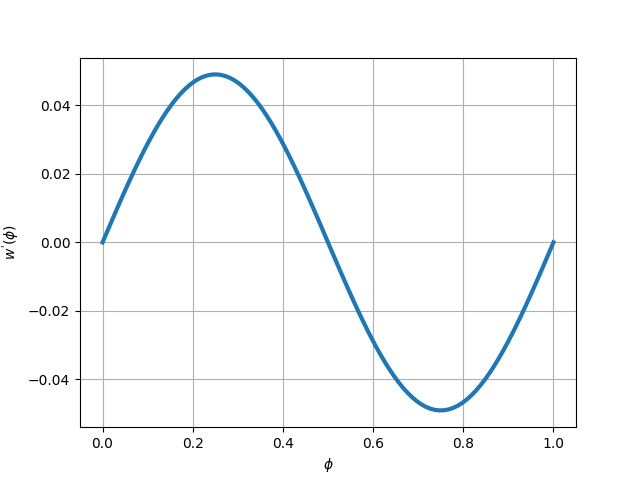}
        \caption{$\frac{dw}{d\phi}$ vs. $\phi$.}
    \end{subfigure}
    \caption{Schematics of the double well potential function $w\left(\phi\right)$ and its derivative $\frac{dw}{d\phi}$ as used in this study.}
    \label{fig:doublewell potential}
\end{figure}

The derivative of $w\left(\phi\right)$ can also be calculated using the automatic differentiation, although the derivation is straightforward in this case.
Schematics of $w\left(\phi\right)$ and $\frac{dw\left(\phi\right)}{d\phi}$ are shown in Figure~\ref{fig:doublewell potential}.
The Cahn-Hilliard model can be readily extended to multi-material topology optimization problems using an array of $\phi_i$ variables, each ranging from 0 to 1 and representing the density of a specific material.
The model described in Equation~\eqref{cahn hilliard model} becomes:
\begin{gather}
    \begin{cases}\label{multi-material cahn hilliard model}
        \frac{\partial \phi_i}{\partial t} = D_i \nabla^2 \mu_i,~~in~\Omega,
        \\
        \\
        \mu_i = \frac{\partial w\left(\phi_i\right)}{\partial \phi_i} - \gamma_i \nabla^2 \phi_i,
        \\
        \\
        \nabla \mu_i = 0,~~in~d\Omega,
        \\
        \\
        \nabla \phi_i = 0,~~in~d\Omega,
    \end{cases}
\end{gather}
where $i=1,2,...,N$ and $N$ is the number of phases. Different coefficients $D_i$, $\gamma_i$ and even $w\left(\phi_i\right)$ can be used for different phases to depict varying phase evolution behaviors, however, here we use same values for all phases for simplicity.
In this study, we update the variables $\phi_i$ in Equation~\eqref{multi-material cahn hilliard model} using the forward Euler method for the time derivative and the finite difference method for the spatial derivatives.
This explicit formulation aligns with the PT method discussed in Section~\ref{the pseudo transient method}, allowing us to fully exploit the parallel capabilities of GPUs. 

\subsection{Workflow of the proposed approach}\label{workflow}

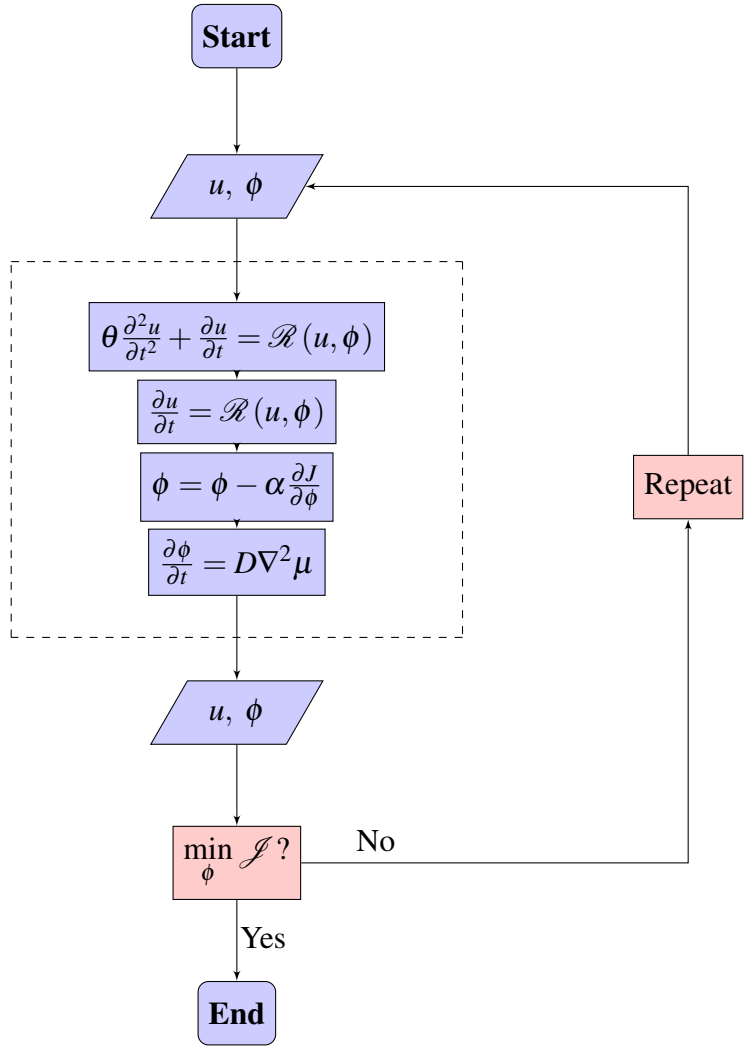
\begin{figure}[hbt!]
\centering
    \begin{tikzpicture}
        \node [terminator, fill=blue!20] at (0,0) (start) {\textbf{Start}};
        \node [data, fill=blue!20] at (0,-2) (theta) {$u,~\phi$};

        \node [process, fill=blue!20] at (0,-4) (PDE_1) {$\theta \frac{\partial^2 u}{\partial t^2} + \frac{\partial u}{\partial t}  = \mathscr{R}\left(u,\phi\right)$};
        \node [process, fill=blue!20] at (0,-5) (PDE_2) {$\frac{\partial u}{\partial t}  = \mathscr{R}\left(u,\phi\right)$};
        \node [process, fill=blue!20] at (0,-6) (sensitivity) {$\phi= \phi- \alpha \frac{\partial J}{\partial \phi}$};
        \node [process, fill=blue!20] at (0,-7) (ch) {$\frac{\partial \phi}{\partial t} = D \nabla^2 \mu$};

        \node [data, fill=blue!20] at (0,-9) (u) {$u,~\phi$};
        \node [process, fill=red!20] at (6,-6) (update) {Repeat};
        \node [process, fill=red!20] at (0,-11) (success) {$\min\limits_{\phi} \mathscr{J} ?$};
        \node [terminator, fill=blue!20] at (0,-13) (end) {\textbf{End}};

        \node[draw=none] at (1.85, -10.7) (no) {No};
        \node[draw=none] at (0.35, -12) (yes) {Yes};
        \path [connector] (start) -- (theta);
        \path [connector] (theta) -- (PDE_1);
        \path [connector] (PDE_1) -- (PDE_2);
        \path [connector] (PDE_2) -- (sensitivity);
        \path [connector] (sensitivity) -- (ch);
        \path [connector] (ch) -- (u);
        \path [connector] (u) -- (success);
        \path [connector] (success) -| (update);
        \path [connector] (update) |- (theta);
        \path [connector] (success) -- (end);

        \draw[dashed] (-3,-3)--(3,-3);
        \draw[dashed] (3,-3)--(3,-8);
        \draw[dashed] (3,-8)--(-3,-8);
        \draw[dashed] (-3,-8)--(-3,-3);

    \end{tikzpicture}
    \caption{A schematic of workflow of our proposed approach. The dashed box outlines one optimization loop, which includes accelerated PT iterations for $u$, PT iterations for $u$, $\phi$ update with sensitivity calculation and $\phi$ update with the Cahn-Hilliard model.}\label{proposed workflow}
\end{figure}

Figure~\ref{proposed workflow} illustrates the workflow of the approach proposed in this study. 
We integrate the methods introduced in the previous sections into a single cohesive algorithm, aiming to solve the topology optimization problem iteratively.
Within each optimization loop, we update the solution $u$ for the governing state equations, update of the design variable $\phi$ using the calculated sensitivities and further refine $\phi$ using the Cahn-Hilliard model.
Therefore, our approach updates the $u$ and $\phi$ simultaneously, whereas traditional workflows update them sequentially.
The optimization process terminates when the value of the objective function, such as thermal compliance, shows a minimal change or when other design specifications have been satisfied.

Note that if we solve the solution $u$ "accurately" within each optimization loop, the approach becomes procedurally identical to the traditional topology optimization approach, as shown in Figure~\ref{a typical workflow of topology optimization}.
This observation naturally introduces a hyperparameter in this study: the balance between the updates of $u$ and $\phi$. Throughout our numerical experiments, we found that the ratio of the iterations for $u$ to the iterations for $\phi$ has an impact on the convergence behavior and the final structures.
More rigorous investigations are needed to fully understand this behavior and this is beyond the scope of current study, thus will not be discussed in detail here.
In the next section, we will solve several topology optimization problems using our proposed approach, aiming to demonstrate its effectiveness and efficiency, while also contributing to the advancement of the topology optimization community.

\FloatBarrier 

\section{Numerical Examples}\label{numerical_examples}

In this section, we present various numerical examples of topology optimization problems to show the effectiveness and efficiency of the proposed approach.
All the examples are implemented using JAX of version 0.4.23, running on Nvidia GeForce RTX 3090 GPUs with CUDA 12.4
and an AMD EPYC 7763 processor. Both GPU and CPU computations are performed using Float32 data type.

\subsection{2D single-material heat conduction optimization}

\begin{figure}[hbt!]
    \centering
    \begin{tikzpicture}
        \draw[color=black] (0.0,0.0) rectangle (4.0,4.0);
        \node at (2.0, 2.0) {$\Omega$};
        \node at (-0.5, 2.0) {$\partial \Omega_D$};
        \node at (2.0, 4.5) {$\partial \Omega_D$};
        \node at (4.5, 2.0) {$\partial \Omega_N$};
        \node at (2.0, -0.5) {$\partial \Omega_N$};
        \node[anchor=north east]
            at (0,0){(0, 0)};
            \node[anchor=south west]
            at (4,4){(4, 4)};

    \end{tikzpicture}
    \caption{Computational domain $\Omega$ for the heat conduction topology optimization problem.
    Dirichlet boundary conditions are applied to the top and left sides $\partial \Omega_D$, while Neumann boundary conditions are applied to the right and bottom sides $\partial \Omega_N$.}
    \label{heat computational domain}
\end{figure}
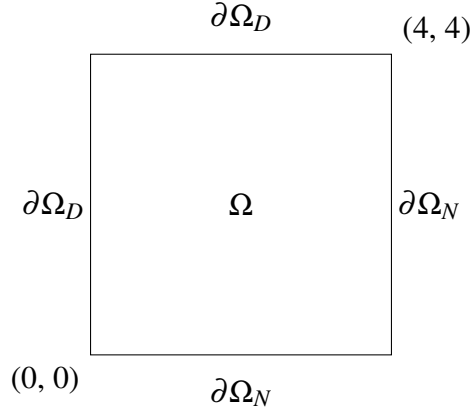

In this section, we present our numerical results for the 2D heat conduction topology optimization problem.
The governing equation is the Poisson's equation discussed in Section~\ref{the pseudo transient method}.
The setup is illustrated in Figure~\ref{heat computational domain}, where Dirichlet boundary conditions with $T=0$ are applied to the top and left sides,  
no-flux Neumann boundary conditions are applied to the right and bottom sides. 
A uniform heat source is present throughout the domain $\Omega$ with $f = 0.01$.
The objective is to distribute a conductive material within the domain to minimize the thermal compliance, as defined in Section~\ref{sensitivity calculation}.
In this case, we consider only one material with a target solid volumetric fraction of $V_s = 0.3$, resulting in a target void volumetric fraction of $V_v = 0.7$.
The conductivity of the material is $\kappa_s = 1.0$. To avoid numerical singularities, the void areas are assigned a conductivity of $\kappa_v = 1 \times 10^{-6}$.
The overall effective heat conductivity can then be calculated as:

\begin{equation}
    \kappa = \kappa_s \phi_s^e + \kappa_v \phi_v^e,\label{overall heat conductivity}
\end{equation}
where $\phi_s$ and $\phi_v$ are respective volumetric variable for the solid and void phases.
The penalty factor $e$ aids in accelerating the separation of phases~\cite{mallon2023neural}. In this study, we use $e = 3$.

The problem is solved on a 512 by 512 structured grid with each side having a length $L = 4.0$. 
The number of degree of freedoms~(DOFs) is $512 \times 512 \times 3$, totaling 786432.
Each optimization loop comprises 500 iterations of the accelerated pseudo-transition solving for the solution $T$, 500 iterations of the PT solving for $T$, one update step of $\phi$s using sensitivity analysis, and one evolution step for $\phi$s using the Cahn-Hilliard model. 
The time step $\Delta t_1$ for the PT iteration is set to $\Delta x^2/4$, where the grid size $\Delta x = 1/512$.
The time step $\Delta t_2$ for the accelerated PT iteration is $\Delta x / 2$ and the time step $\Delta t_3$ for the Cahn-Hilliard model is $500 \Delta x^4$. 
For this and all other examples, $\gamma = 3 \times 10^{-5}$ is used in the Cahn-Hilliard model as described in Equation~\eqref{cahn hilliard model}.

In the overall loss function $J$ formulation, we assign a weight of $\alpha_v = 10^5$ to the volumetric loss component $J_v$, and a weight of $\alpha_1 = 10^4$ to the partition of unity loss component $J_1$. 
For the thermal compliance loss component, we employ a normalized format, 
\begin{equation}
    J^i_h = \alpha_h \frac{\frac{\partial J_h}{\partial \phi^n_i}}{\max\left(\left|\frac{\partial J_h}{\partial \phi_i^n}\right|\right)},
\end{equation}
where $\alpha_h$ is set to 0.1. The update for the density variables using the sensitivity can then be written as:
\begin{equation}
    \phi_i^{n+1} = \phi_i^{n} - \left(J^i_h + \alpha_v \frac{\partial J_v}{\partial \phi_i^n} + \alpha_1 \frac{\partial J_1}{\partial \phi_i^n}\right),
\end{equation}


\FloatBarrier 

\begin{figure}[hbt!]
    \centering
    \begin{subfigure}{0.4\textwidth}
        \centering
        \includegraphics[width=0.9\linewidth]{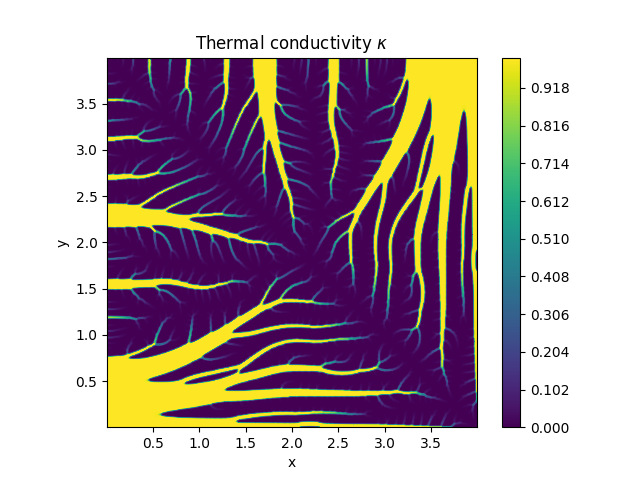}
        \caption{Thermal conductivity $\kappa$.}
    \end{subfigure}
    \begin{subfigure}{0.4\textwidth}
        \centering
        \includegraphics[width=0.9\linewidth]{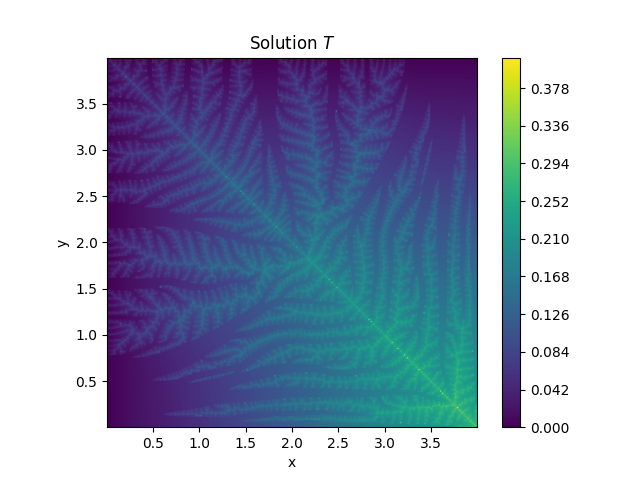}
        \caption{Solution $T$.}
    \end{subfigure}
    \caption{The optimized distribution of $\kappa$ after 5000 iterations and the corresponding solution $T$.}
    \label{fig:thermal solutions}
\end{figure}

\begin{figure}[hbt!]
    \centering
    \begin{subfigure}{0.4\textwidth}
        \centering
        \includegraphics[width=0.9\linewidth]{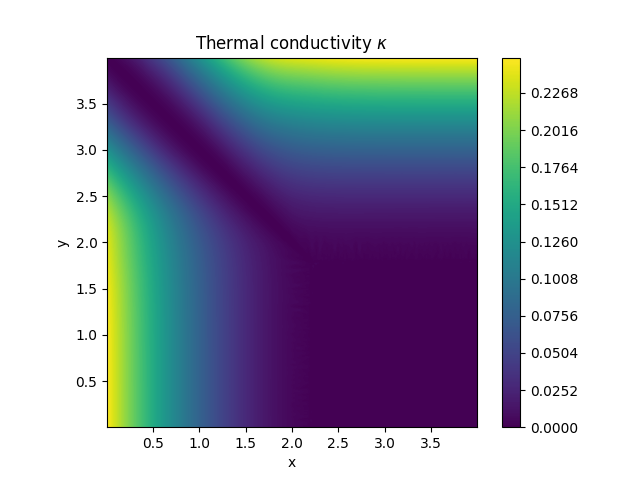}
        \caption{$\kappa$ after 100 iterations.}
    \end{subfigure}
    \begin{subfigure}{0.4\textwidth}
        \centering
        \includegraphics[width=0.9\linewidth]{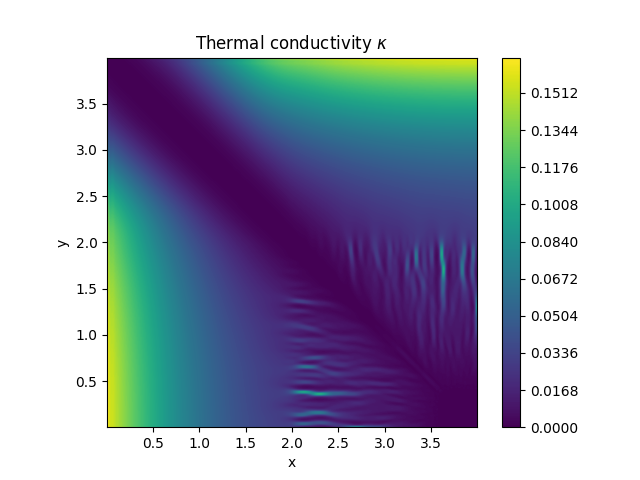}
        \caption{$\kappa$ after 500 iterations.}
    \end{subfigure}
    \begin{subfigure}{0.4\textwidth}
        \centering
        \includegraphics[width=0.9\linewidth]{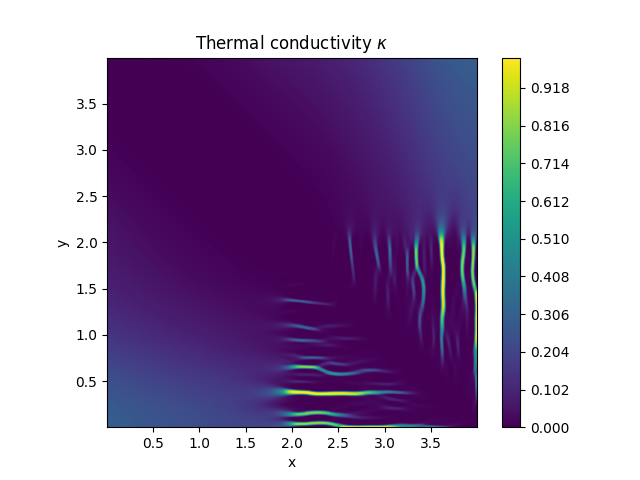}
        \caption{$\kappa$ after 700 iterations.}
    \end{subfigure}
    \begin{subfigure}{0.4\textwidth}
        \centering
        \includegraphics[width=0.9\linewidth]{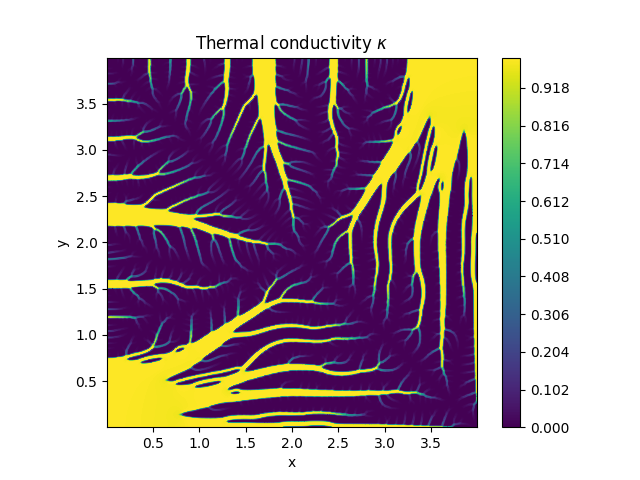}
        \caption{$\kappa$ after 2000 iterations.}
    \end{subfigure}
    \caption{The distribution of $\kappa$ after different number of iterations: 100, 500, 700, 2000 iterations.}
    \label{fig:different thermal solutions}
\end{figure}

\begin{figure}[hbt!]
    \centering
    \begin{subfigure}{0.4\textwidth}
        \centering
        \includegraphics[width=0.9\linewidth]{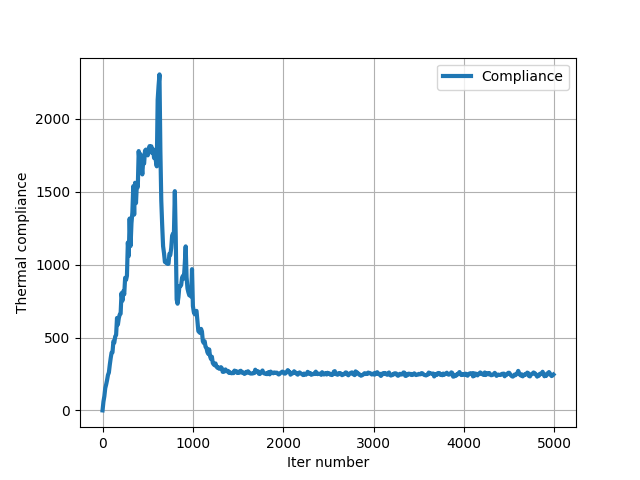}
        \caption{Thermal compliance.}
    \end{subfigure}
    \begin{subfigure}{0.4\textwidth}
        \centering
        \includegraphics[width=0.9\linewidth]{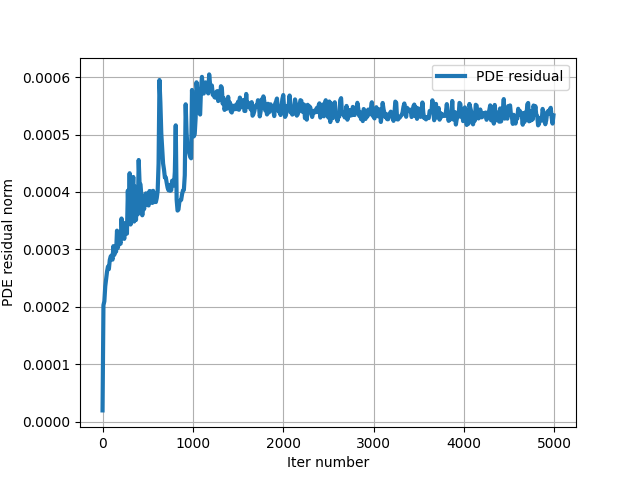}
        \caption{The governing PDE residual.}
    \end{subfigure}
    \begin{subfigure}{0.4\textwidth}
        \centering
        \includegraphics[width=0.9\linewidth]{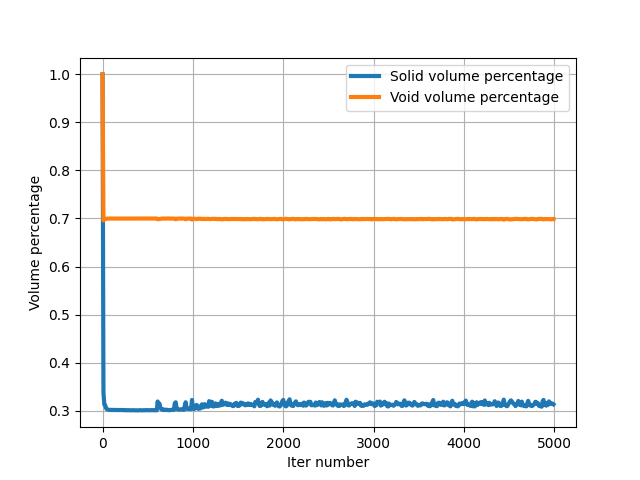}
        \caption{Volumetric fraction of phases.}
    \end{subfigure}
    \caption{Convergence plots showing the thermal compliance, PDE residual and volumetric fractions of different phases throughout the iteration process.}
    \label{fig:thermal convergence plots}
\end{figure}

We begin with initial conditions of $T = 0$ throughout the domain and $\phi_i = 1.0$ for both phases. Figure~\ref{fig:thermal solutions} illustrates the optimized thermal conductivity $\kappa$ and its corresponding solution $T$.
The two phases (solid and void) have separated after 5000 optimization iterations.
Figure~\ref{fig:different thermal solutions} shows the evolution of $\kappa$, with the distribution noticeably separated after 2000 iterations.
Here we define the residual norm corresponding to the Equation~\eqref{Poisson equation} as:
\begin{align}
    r_{i,j} = \nabla \cdot \left(\kappa \nabla T\right)_{i,j} + f_{i,j},
    \\
    r_{PDE} = \frac{\left(\sum_{i,j} r_{i,j}^2\right)^{1/2}}{\mathcal N}, 
\end{align}
where $r_{i,j}$ is the residual of each grid and $\mathcal N$ is the total number of grids.
Since we have converted the topology optimization problem into solving a set of time dependent PDEs problems, the residual of the Poisson equation converges as the iteration progresses shown in Figure~\ref{fig:thermal convergence plots}.
Moreover, the thermal compliance stabilizes after an initial spike, while adhering to the imposed volumetric constraints.
Completing 5000 iterations takes 51 seconds on GPUs and 2173 seconds on CPUs.
A significant speedup is observed when using the same code on GPUs compared to CPUs.

\FloatBarrier 

\FloatBarrier 

\FloatBarrier 

\begin{figure}[hbt!]
    \centering
    \begin{subfigure}{0.4\textwidth}
        \centering
        \includegraphics[width=0.9\linewidth]{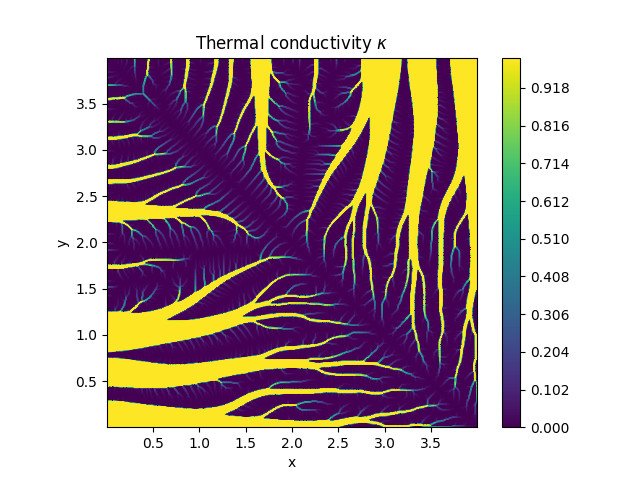}
        \caption{$\Delta t_3 = 100 \Delta x^4$.}
    \end{subfigure}
    \begin{subfigure}{0.4\textwidth}
        \centering
        \includegraphics[width=0.9\linewidth]{figures/2d_k_sol.png}
        \caption{$\Delta t_3 = 500 \Delta x^4$.}
    \end{subfigure}
    \begin{subfigure}{0.4\textwidth}
        \centering
        \includegraphics[width=0.9\linewidth]{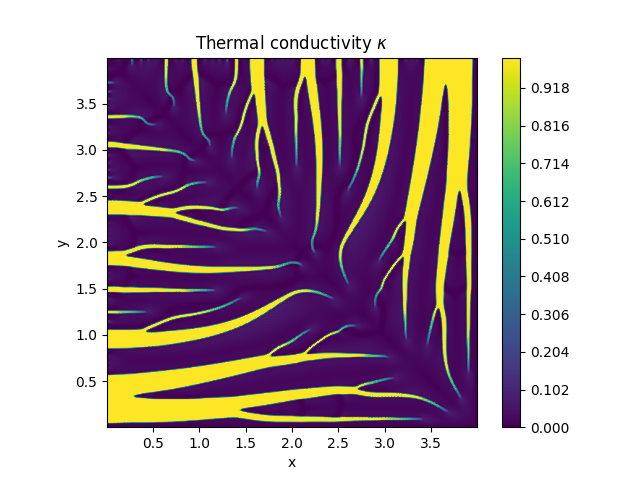}
        \caption{$\Delta t_3 = 1000 \Delta x^4$.}
    \end{subfigure}
    \caption{Distribution of $\kappa$ after 5000 iterations using different time steps $\Delta t_3$ in the Cahn-Hilliard model.}
    \label{fig:k result with different alpha values}
\end{figure}

Figure~\ref{fig:k result with different alpha values} illustrates the final distributions of $\kappa$ using different $\Delta t_3$ in the Cahn-Hilliard model after 5000 iterations.
A larger $\Delta t_3$ results in a smoother distribution with fewer branches. This is an active research area in topology optimization~\cite{huang2024derivable},
focusing on effective size control of distributed materials to incorporate additional manufacturing consideration into the topology optimization process. 
In this work, we present some empirical observations and leave a more thorough analysis for further research. 

\FloatBarrier 

\subsection{2D multi-material MBB beam}

In this section, we present our numerical results for the 2D Messerschmitt-Bolkow-Blohm~(MBB) beam topology optimization.
The governing equation is the classical elasticity equation in Equation~\eqref{momentum balance}.

\begin{gather}
    \begin{cases}\label{momentum balance}
        \frac{\sigma_{xx}}{\partial x} + \frac{\sigma_{xy}}{\partial y} = f_x,
        \\
        \\
        \frac{\sigma_{yx}}{\partial x} + \frac{\sigma_{yy}}{\partial y} = f_y.
    \end{cases}
\end{gather}

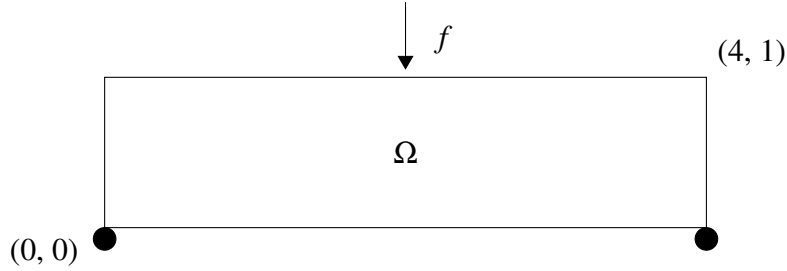
\begin{figure}[hbt!]
    \centering
    \begin{tikzpicture}
        \draw[color=black] (0.0,0.0) rectangle (8.0,2.0);
        \node at (4.0, 1.0) {$\Omega$};
        \draw[color=black, fill=black] (0.0,-0.15) circle (0.15);
        \draw[color=black, fill=black] (8.0,-0.15) circle (0.15);
        \draw[->, color=black, -triangle 60] (4.0, 3.0) -- (4.0, 2.1);
        \node at (4.5, 2.5) {$f$};
        \node[anchor=north east]
            at (-0.2,0){(0, 0)};
            \node[anchor=south west]
            at (8,2){(4, 1)};
    \end{tikzpicture}
    \caption{Computational domain $\Omega$ for the MBB beam topology optimization problem.}
    \label{2d mbb computational domain}
\end{figure}

Figure~\ref{2d mbb computational domain} illustrates the setup as a $4 \times 1$ computational domain. 
Both bottom corners rest on rollers, preventing vertical displacements. A downward force $f = 1.0$ is exerted at the center of the top side  
while the remaining sides have free surface boundary conditions. 
The objective is to distribute five materials, each with Young's modulus values $E_i$ of $[1.0, 0.775, 0.55, 0.325, 0.1]$, within the domain $\Omega$ to minimize the mechanical compliance $J_m$, defined as:

\begin{equation}\label{thermal compliance definition}
    J_m = \int_{\Omega} \sigma : \epsilon d\Omega,
\end{equation}
where $\sigma$ is the stress tensor and $\epsilon$ is the strain tensor.
Poisson's ratio $\gamma = 0.3$ is used for all materials. The void area is consider to have Young's modulus of $1 \times 10^{-6}$ to avoid numerical singularities.
The target volumetric fraction $V_i$ of each material is set to be 0.08, leaving the void volumetric fraction $V_v$ at 0.6. 
The overall Young's modulus can then be calculated as:

\begin{equation}\label{overall Young's modulus}
    E = \sum_i E_i \phi_i^e,~~i=1,2,3,4,5,
\end{equation}
where the penalty parameter $e = 3$ is used here.

\FloatBarrier 

In this work, we discretize the elasticity equation using the finite difference method and update the solution using both the PT and accelerated PT method.
$\sigma_{xx}$, $\sigma_{xy}$, $\sigma_{yx}$ and $\sigma_{yy}$ in Equation~\eqref{momentum balance} can be further expressed as:

\begin{gather}\label{stress definition}
    \begin{cases}
        \sigma_{xx} = \left(\lambda + 2 \mu \right) \frac{\partial u_x}{\partial x} + \lambda \frac{\partial u_y}{\partial y},
        \\
        \\
        \sigma_{yy} = \lambda \frac{\partial u_x}{\partial x} + \left(\lambda + 2 \mu\right) \frac{\partial u_y}{\partial y},
        \\
        \\
        \sigma_{xy} = \sigma_{yx} = \mu \left(\frac{\partial u_x}{\partial y} + \frac{\partial u_y}{\partial x}\right),
    \end{cases}
\end{gather}
$u_x$ and $u_y$ are displacements in $x$ and $y$ directions. The Lame parameters $\lambda$ and $\mu$ can be calculated with the overall Young's modulus $E$ and Poisson's ratio $\gamma$ as:

\begin{gather}
    \begin{cases}
        \lambda = \frac{E \gamma}{\left(1 + \gamma\right) \left(1 - 2 \gamma\right)},
        \\
        \\
        \mu = \frac{E}{2 \left(1 + \gamma\right)}.
    \end{cases}
\end{gather}
The PT and accelerated PT method can be similarly formulated for both variables $u_x$ and $u_y$:

\begin{gather}
    \begin{cases}
        \frac{u_x^{n+1} - u_x^{n}}{\Delta t_1} = \frac{\sigma_{xx}^{n}}{\partial x} + \frac{\sigma_{xy}^{n}}{\partial y} - f_x,
        \\
        \\
        \frac{u_y^{n+1} - u_y^{n}}{\Delta t_1} = \frac{\sigma_{yx}^{n}}{\partial x} + \frac{\sigma_{yy}^{n}}{\partial y} - f_y,
    \end{cases}
\end{gather}

\begin{gather}
    \begin{cases}
        \theta \frac{u_x^{n+1} - 2u_x^{n} + u_x^{n-1}}{\Delta t_2^2} + \frac{u_x^{n+1} - u_x^{n}}{\Delta t_2} = \frac{\sigma_{xx}^{n}}{\partial x} + \frac{\sigma_{xy}^{n}}{\partial y} - f_x,
        \\
        \\
        \theta \frac{u_y^{n+1} - 2u_y^{n} + u_y^{n-1}}{\Delta t_2^2} + \frac{u_y^{n+1} - u_y^{n}}{\Delta t_2} = \frac{\sigma_{yx}^{n}}{\partial x} + \frac{\sigma_{yy}^{n}}{\partial y} - f_y,
    \end{cases}
\end{gather}
here we use $\theta = 1$, with $\Delta t_1 = \Delta x^2 / 4$ and $\Delta t_2 = \Delta x / 2$.
The problem is solved on a $513 \times 128$ uniform mesh, resulting in a total of 525312 DOFs. 
The initial conditions for displacements are $u_x = 0$ and $u_y = 0$ throughout the domain, and an initial value of 0.5 is used for all $\phi$s values.
Each optimization loop comprises 20 iterations of the accelerated pseudo-transition solving for $u_x$ and $u_y$, 20 iterations of the PT solving for $u_x$ and $u_y$, one update step of $\phi$s using sensitivity analysis, and one evolution step for $\phi$s using the multi-material Cahn-Hilliard model in Equation~\eqref{multi-material cahn hilliard model} with $\Delta t_3 = 500\Delta x^4$. 
The sensitivity calculation is similar with the heat conduction example, with slightly different weights: $\alpha_v = 10^4$ and $\alpha_1 = 10^3$. 

\begin{figure}[hbt!]
    \centering
    \begin{subfigure}{0.45\textwidth}
        \centering
        \includegraphics[width=1.0\linewidth]{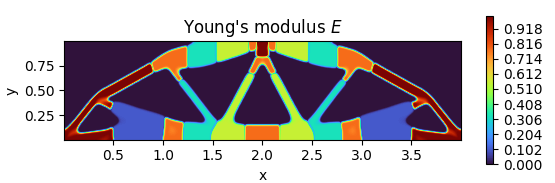}
        \caption{Distribution of $E$.}
    \end{subfigure}
    \begin{subfigure}{0.45\textwidth}
        \centering
        \includegraphics[width=1.0\linewidth]{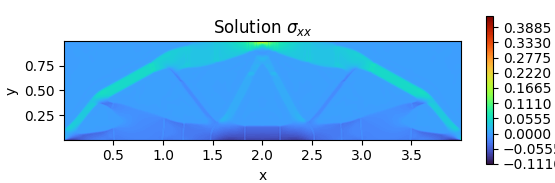}
        \caption{Distribution of $\sigma_{xx}$.}
    \end{subfigure}
    \begin{subfigure}{0.45\textwidth}
        \centering
        \includegraphics[width=1.0\linewidth]{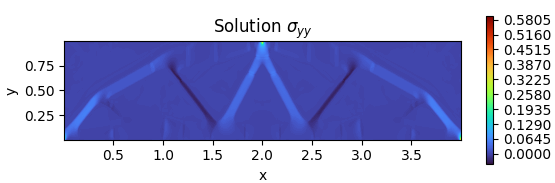}
        \caption{Distribution of $\sigma_{yy}$.}
    \end{subfigure}
    \begin{subfigure}{0.45\textwidth}
        \centering
        \includegraphics[width=1.0\linewidth]{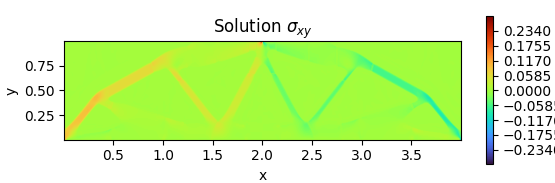}
        \caption{Distribution of $\sigma_{xy}$.}
    \end{subfigure}
    \caption{Overall Young's modulus and stresses after \num{200000} iterations.}
    \label{fig:E results and stresses}
\end{figure}

\FloatBarrier 

\begin{figure}[hbt!]
    \centering
    \begin{subfigure}{0.45\textwidth}
        \centering
        \includegraphics[width=1.0\linewidth]{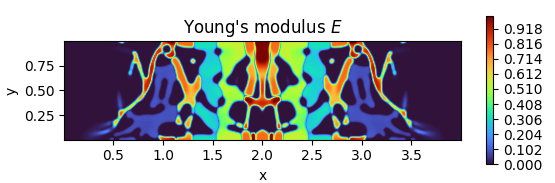}
        \caption{$E$ after \num{25000} iterations.}
    \end{subfigure}
    \begin{subfigure}{0.45\textwidth}
        \centering
        \includegraphics[width=1.0\linewidth]{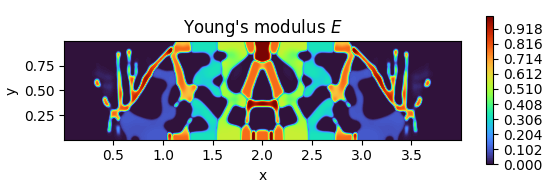}
        \caption{$E$ after \num{50000} iterations.}
    \end{subfigure}
    \begin{subfigure}{0.45\textwidth}
        \centering
        \includegraphics[width=1.0\linewidth]{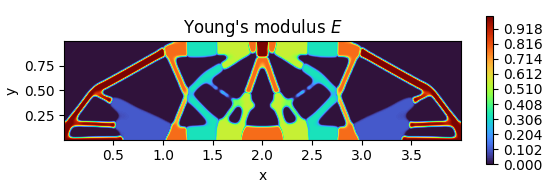}
        \caption{$E$ after \num{75000} iterations.}
    \end{subfigure}
    \begin{subfigure}{0.45\textwidth}
        \centering
        \includegraphics[width=1.0\linewidth]{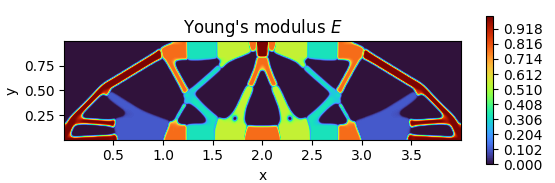}
        \caption{$E$ after \num{100000} iterations.}
    \end{subfigure}
    \caption{The evolution of the overall Young's modulus.}
    \label{fig:different mechanical solutions}
\end{figure}

\FloatBarrier 

Figure~\ref{fig:E results and stresses} displays the overall Young's modulus $E$ across the domain after \num{200000} iterations, along with the corresponding stress distributions. 
Different materials are segregated with values approaching either 0 or 1 across the domain, allowing us to use the distribution of $E$ to illustrate the allocation of each material.
The stresses exhibit a smooth distribution throughout the solid structures, as depicted in Figure~\ref{fig:E results and stresses} (b-d).
In addition, it is evident that harder materials are concentrated in regions of high stress, such as at the load point and the bottom two corners. This distribution highlights how our approach contributes to achieving more effective and efficient designs with different materials.

\begin{figure}[hbt!]
    \centering
    \begin{subfigure}{0.4\textwidth}
        \centering
        \includegraphics[width=0.9\linewidth]{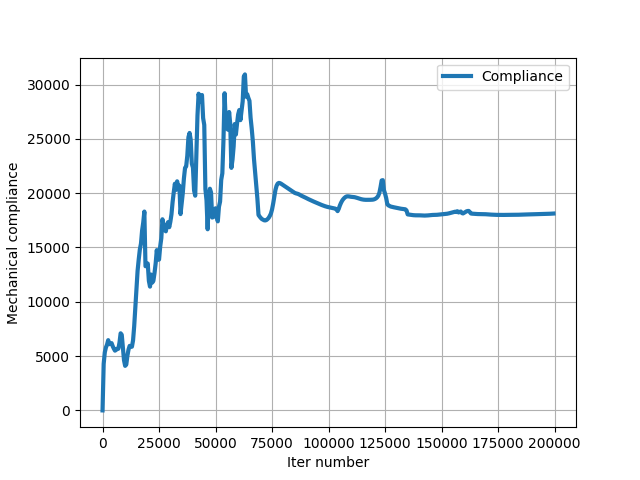}
        \caption{Mechanical compliance.}
    \end{subfigure}
    \begin{subfigure}{0.4\textwidth}
        \centering
        \includegraphics[width=0.9\linewidth]{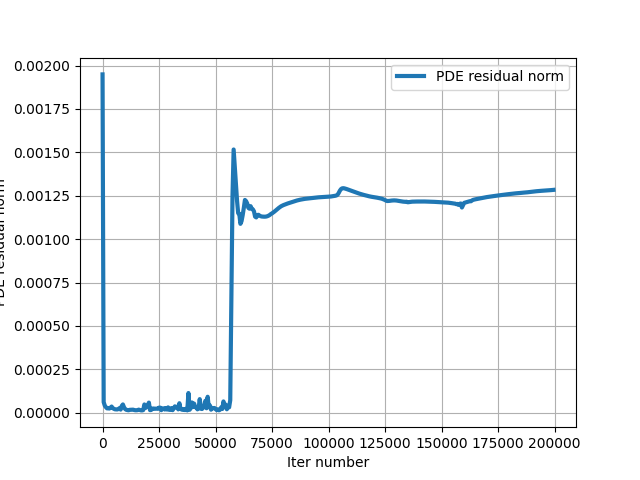}
        \caption{The governing PDE residual.}
    \end{subfigure}
    \begin{subfigure}{0.4\textwidth}
        \centering
        \includegraphics[width=0.9\linewidth]{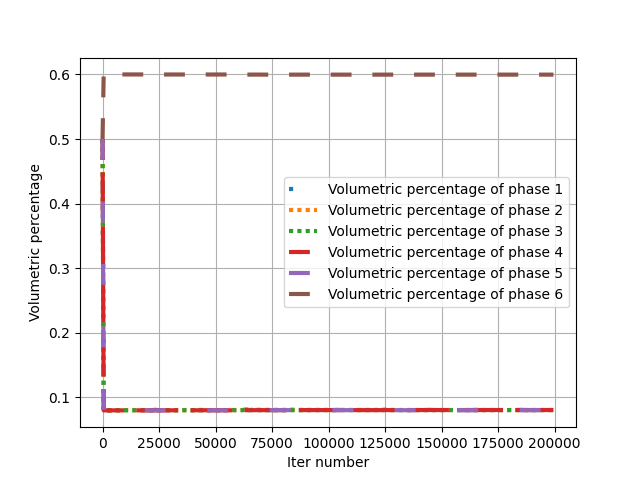}
        \caption{Volumetric fraction of phases.}
    \end{subfigure}
    \caption{Convergence plots showing the mechanical compliance, PDE residual and volumetric fractions of different phases throughout the iteration process.}
    \label{fig:mechanical convergence plots}
\end{figure}

Figure~\ref{fig:different mechanical solutions} illustrates the evolution process for $E$ of the overall structure. Through the use of a relatively small number of PT and accelerated PT iterations in each optimization loop, it vividly demonstrates how materials are integrated and flow into their final allocations.
Figure~\ref{fig:mechanical convergence plots} illustrates that the final solution converges after approximately \num{150000} iterations, ensuring that all materials conform to the specified volumetric fractions. 
Completing \num{200000} iterations takes 176 seconds on GPUs and 8774 seconds on CPUs.
Again, this results in an approximately 50-fold acceleration simply by running the same code on GPUs.

\FloatBarrier 

\subsection{3D multi-material cantilever beam}

\begin{figure}[hbt!]
    \centering
    \includegraphics[width=0.5\linewidth]{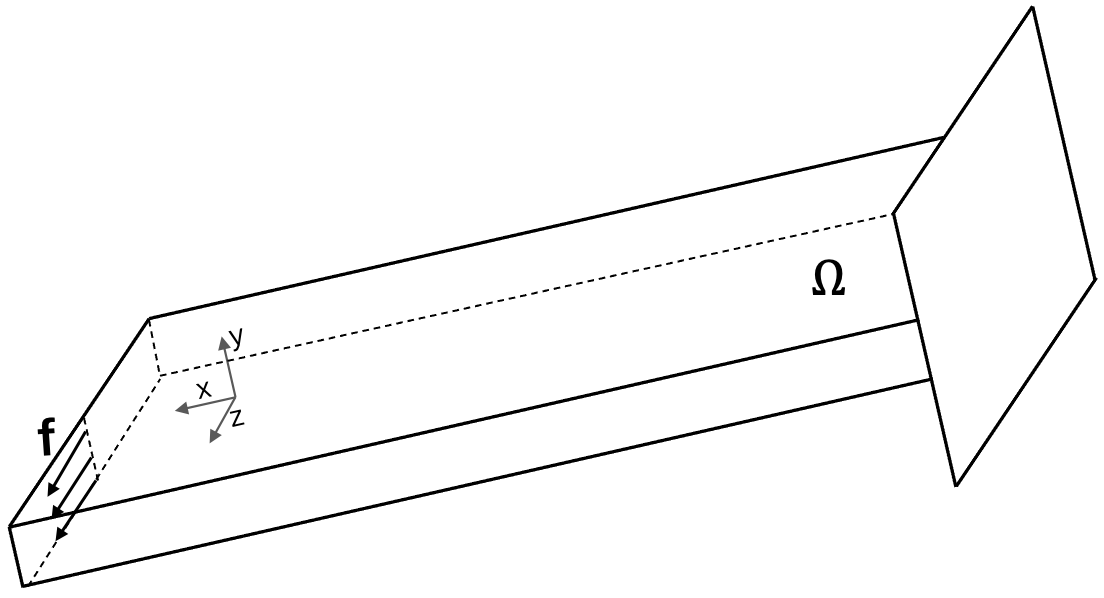}
    \caption{Computational domain $\Omega$ for the 3D cantilever beam problem.}
    \label{fig:3d cantilever domain}
\end{figure}

This section presents our results for the topology optimization of a 3D multi-material cantilever beam.
A schematic is shown in Figure~\ref{fig:3d cantilever domain} with lengths 2, 2/15 and 2/3 in x, y and z directions, respectively. The problem is solved on a $256 \times 17 \times 85$ structured grid.
A uniform force $f = 1$, in the positive z-direction, is applied to the centerline of the left side of the beam, with the right side fixed and the remaining surfaces are left free.
We aim to distribute three materials with $E = 1.0, 0.6$ and 0.2 within the domain to minimize the mechanical compliance, subject to a volumetric fraction 0.1 for each material.
This results in the volumetric fraction of 0.7 for the void volume. To avoid numerical singularities, $E = 1 \times 10^{-6}$ is used for the void area.

\FloatBarrier 

\begin{figure}[hbt!]
    \centering
    \begin{subfigure}{0.4\textwidth}
        \centering
        \includegraphics[width=0.9\linewidth]{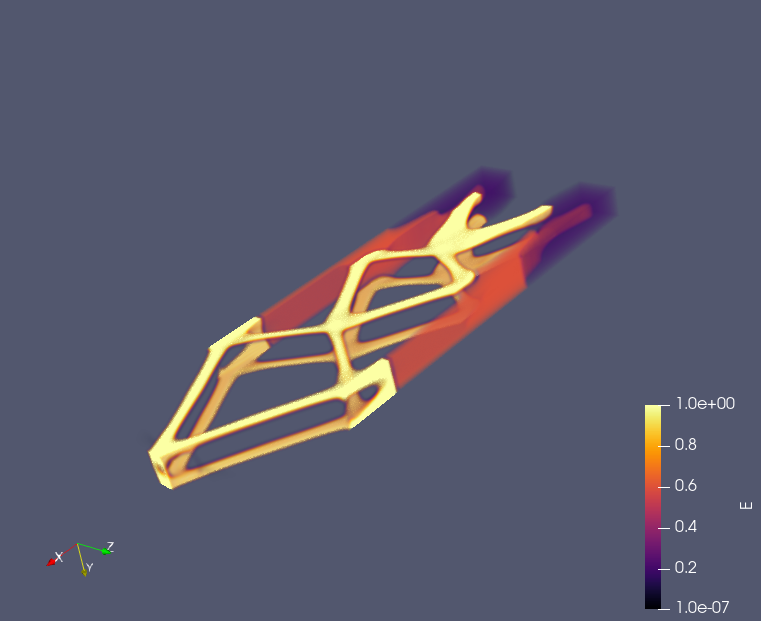}
        \caption{View 1.}
    \end{subfigure}
    \begin{subfigure}{0.4\textwidth}
        \centering
        \includegraphics[width=0.9\linewidth]{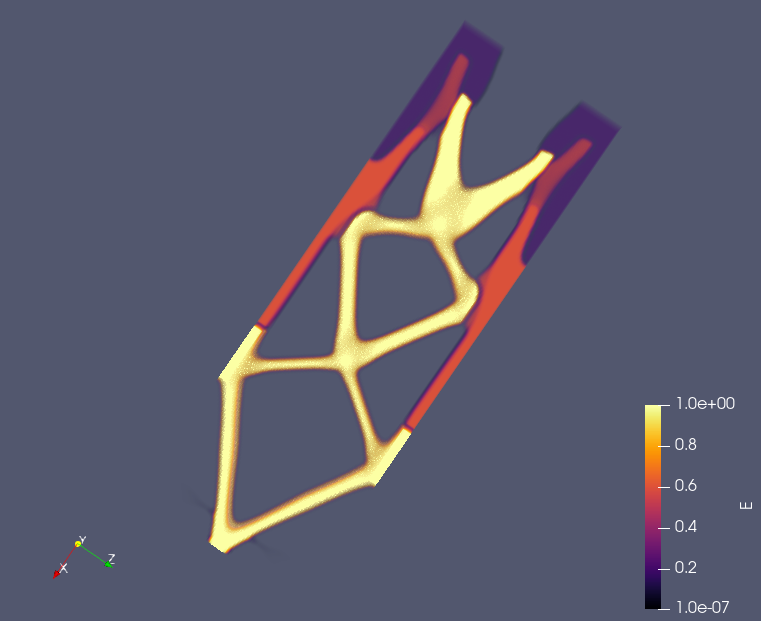}
        \caption{View 2.}
    \end{subfigure}
    \begin{subfigure}{0.4\textwidth}
        \centering
        \includegraphics[width=0.9\linewidth]{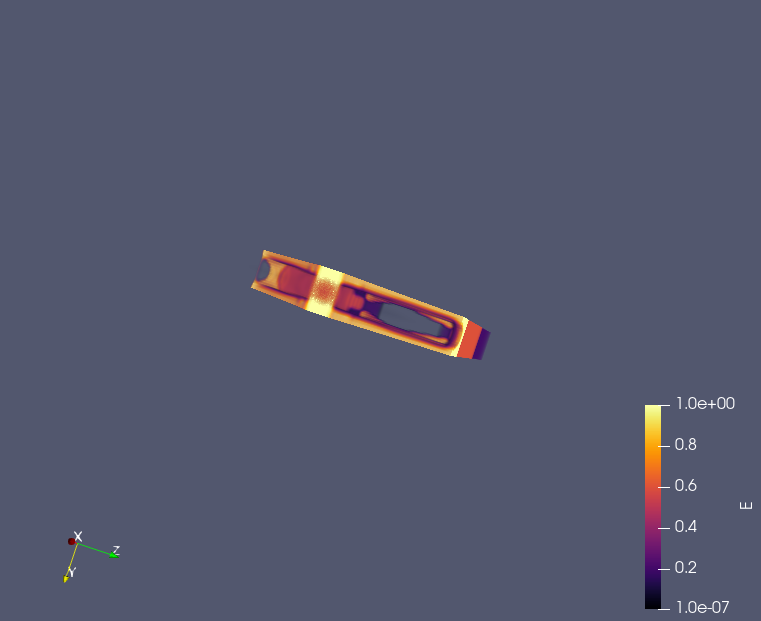}
        \caption{View 3.}
    \end{subfigure}
    \begin{subfigure}{0.4\textwidth}
        \centering
        \includegraphics[width=0.9\linewidth]{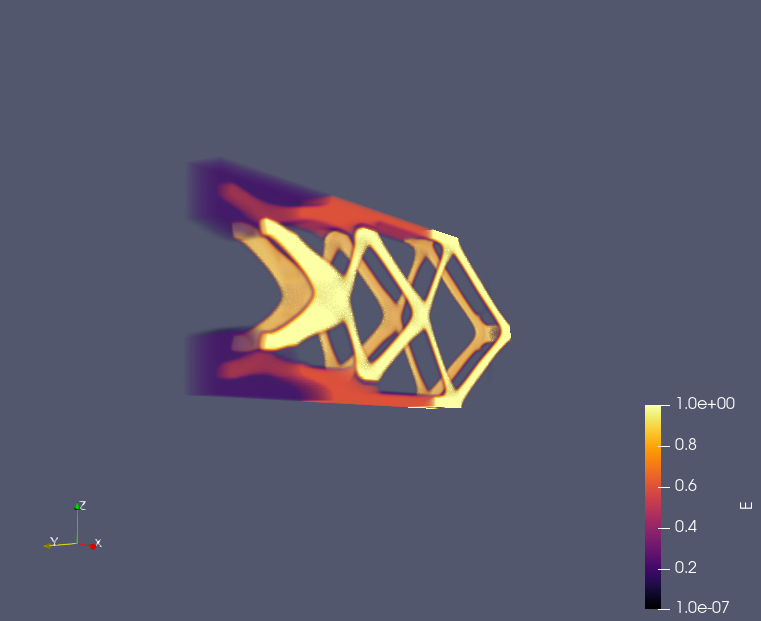}
        \caption{View 4.}
    \end{subfigure}
    \caption{Four views of the optimized 3D cantilever beam.}
    \label{fig:3d cantilever views}
\end{figure}

The workflow and implementation are essentially the same as the 2D MBB problem except this time the linear elasticity equation is solved for three displacements $u_x$, $u_y$ and $u_z$. 
Therefore, the DOFs for this problem amount to $\num{369920} \times 6= \num{2219520}$.
For this problem, each optimization loop comprises 100 iterations of the accelerated pseudo-transition iterations, 100 iterations of the PT iterations, one update step of $\phi$s using sensitivity analysis, and one evolution step for $\phi$s using the Cahn-Hilliard model. 
In the PT iteration, we use $\Delta t_1 = \Delta x^2 / 6$ while in the accelerated PT iteration, $\Delta t_2 = \Delta x / 2$ is employed. Additionally, $\Delta t_3 = 500\Delta x^4$ is used in the Cahn-Hilliard model. 
The weights used in the sensitivity calculation are identical to those used in the 2D MBB problem.

Four views depicting the distribution of $E$ after 4000 optimization steps are shown in Figure~\ref{fig:3d cantilever views}. 
It demonstrates how the materials are distributed to form a resilient structure while adhering to volumetric constraints. 
The entire optimization process takes 796 seconds on GPUs, compared to approximately 36800 seconds on CPUs.

\FloatBarrier 

\subsection{3D single-material drone skeleton}

\begin{figure}[hbt!]
    \centering
    \includegraphics[width=0.8\linewidth]{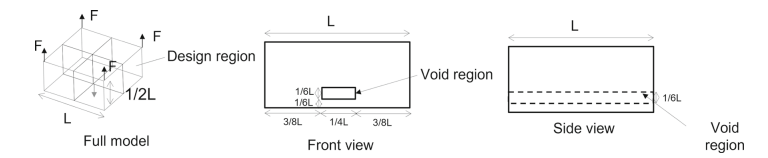}
    \caption{Computational domain for the drone skeleton design from paper by Li et al.\cite{li2023convolution}.}
    \label{fig:3d drone domain}
\end{figure}

This section presents our results for topology optimization of a drone skeleton using a single solid material. 
A schematic of the design specification is shown in Figure~\ref{fig:3d drone domain}, detailed in the paper by Li et al.\cite{li2023convolution}. 
In short, the four top corners are constrained against rollers to prevent vertical displacements, while a point load is applied at the center of the bottom side. All other surfaces are left free.
One noteworthy aspect here is the additional volumetric constraint: that there is a fraction of the domain must remain void to accommodate accessories such as batteries. 
This constraint can be easily incorporated by modifying the sensitivity definition:

\begin{equation}\label{drone sensitivity}
    \phi_i^{n+1} = \phi_i^{n} - \left(J^i_m + \alpha_v \frac{\partial J_v}{\partial \phi_i^n} + \alpha_1 \frac{\partial J_1}{\partial \phi_i^n} + \alpha_b \frac{\partial J_b}{\partial \phi_i^n}\right),
\end{equation}
where $J_b$ similar to $J_v$, represents the volumetric constraint for the void region, defined as:

\begin{equation}
    J_b = \sum_i \left(\int_{\Omega_b} \phi_i d\Omega_b - V^b_i\right)^2,
\end{equation}
where $\Omega_b$ is the void region, with $V^b_1 = 0.0$ indicating solid material and $V^b_2 = 1.0$ indicating void volume.

\FloatBarrier 

\begin{figure}[hbt!]
    \centering
    \begin{subfigure}{0.4\textwidth}
        \centering
        \includegraphics[width=0.9\linewidth]{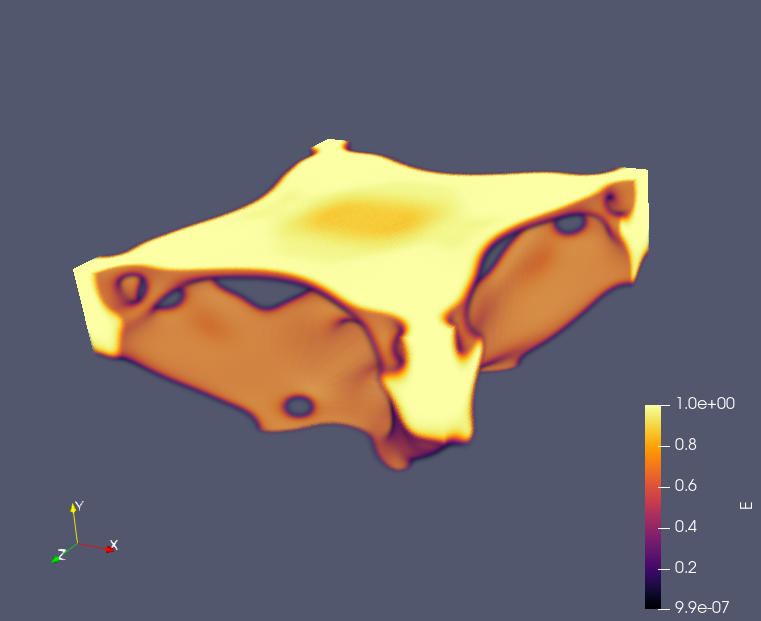}
        \caption{View 1.}
    \end{subfigure}
    \begin{subfigure}{0.4\textwidth}
        \centering
        \includegraphics[width=0.9\linewidth]{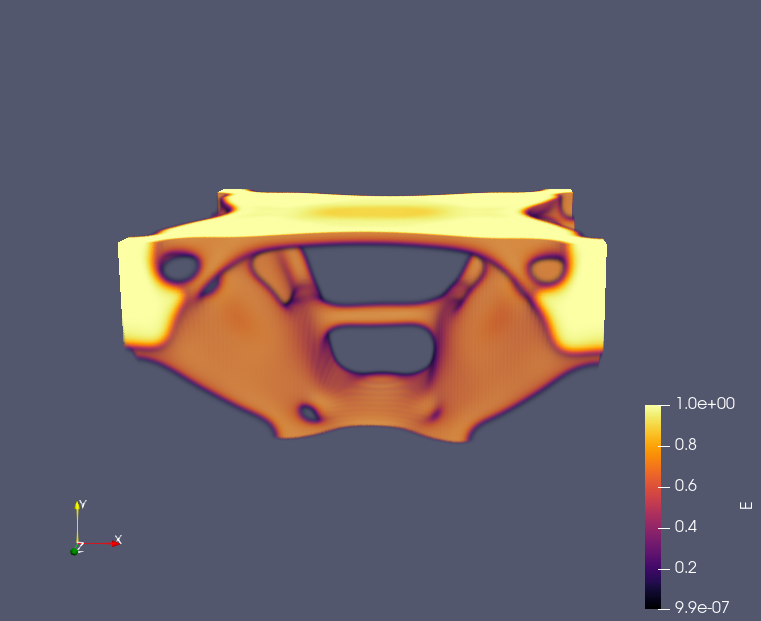}
        \caption{View 2.}
    \end{subfigure}
    \begin{subfigure}{0.4\textwidth}
        \centering
        \includegraphics[width=0.9\linewidth]{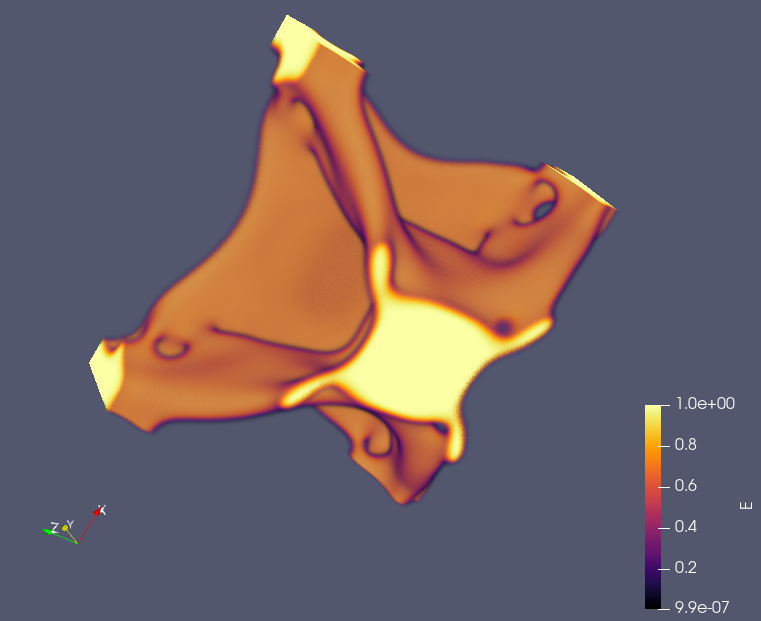}
        \caption{View 3.}
    \end{subfigure}
    \begin{subfigure}{0.4\textwidth}
        \centering
        \includegraphics[width=0.9\linewidth]{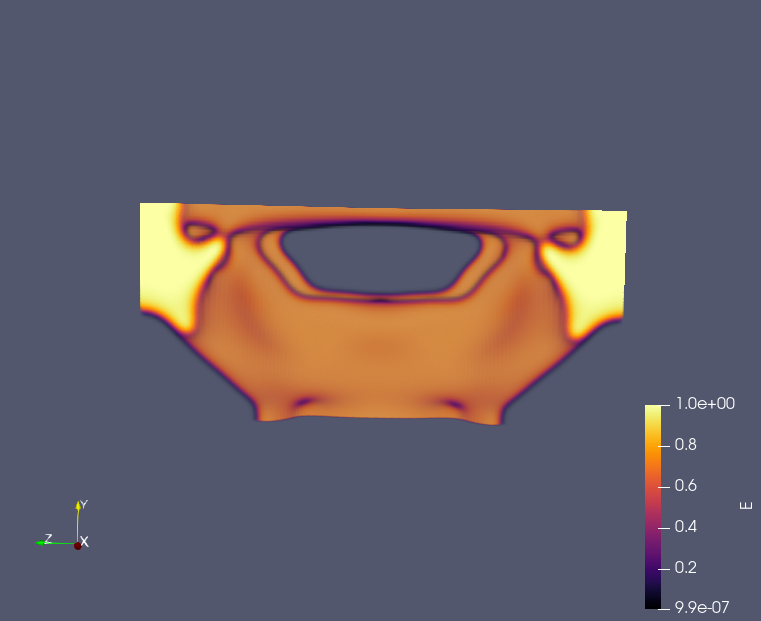}
        \caption{View 4.}
    \end{subfigure}
    \caption{Four views of the drone structure with 0.2 void volumetric fraction.}
    \label{fig:drone 020 fraction views}
\end{figure}

\begin{figure}[hbt!]
    \centering
    \begin{subfigure}{0.4\textwidth}
        \centering
        \includegraphics[width=0.9\linewidth]{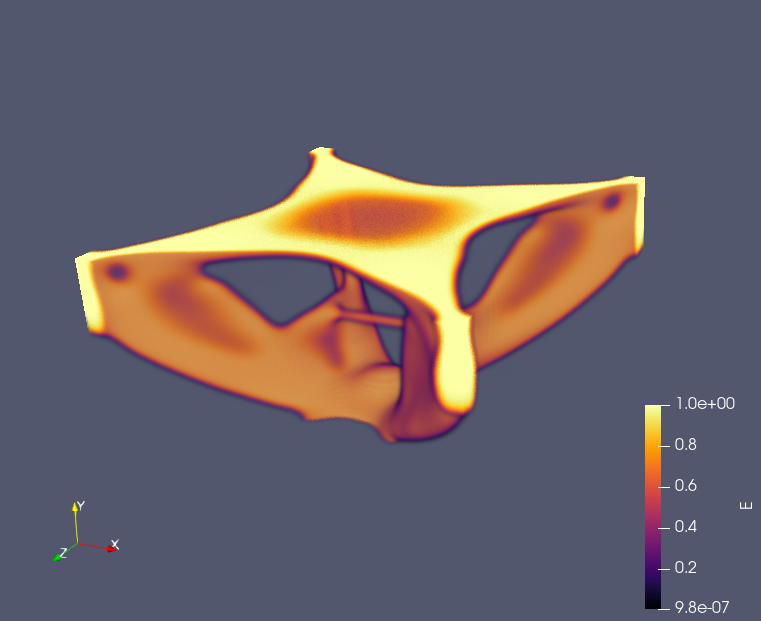}
        \caption{View 1.}
    \end{subfigure}
    \begin{subfigure}{0.4\textwidth}
        \centering
        \includegraphics[width=0.9\linewidth]{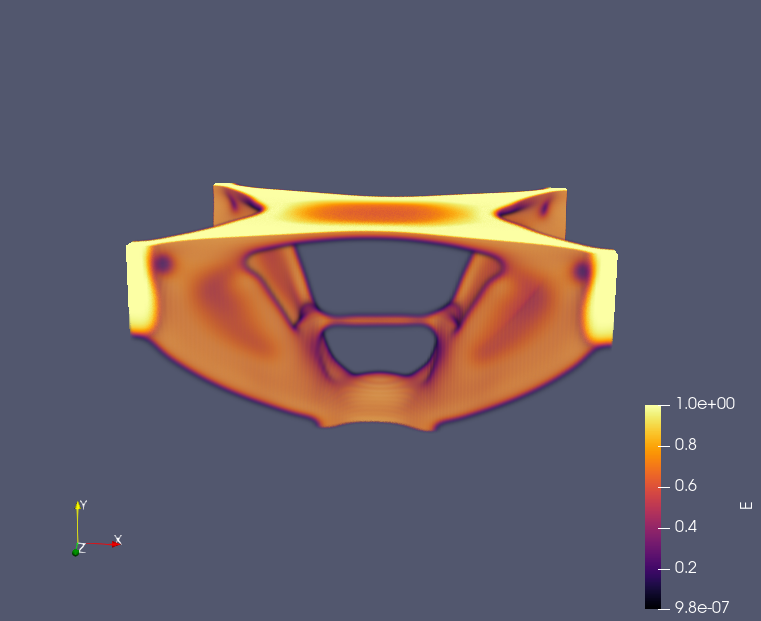}
        \caption{View 2.}
    \end{subfigure}
    \begin{subfigure}{0.4\textwidth}
        \centering
        \includegraphics[width=0.9\linewidth]{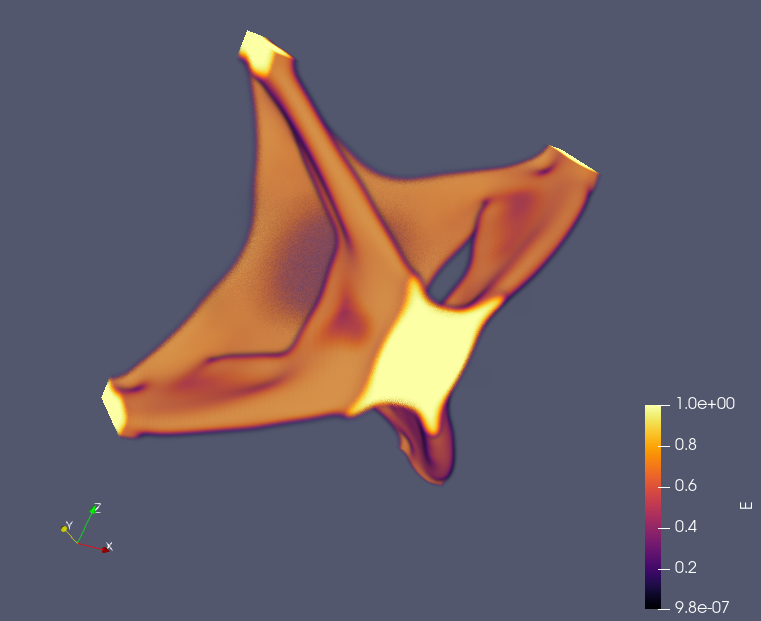}
        \caption{View 3.}
    \end{subfigure}
    \begin{subfigure}{0.4\textwidth}
        \centering
        \includegraphics[width=0.9\linewidth]{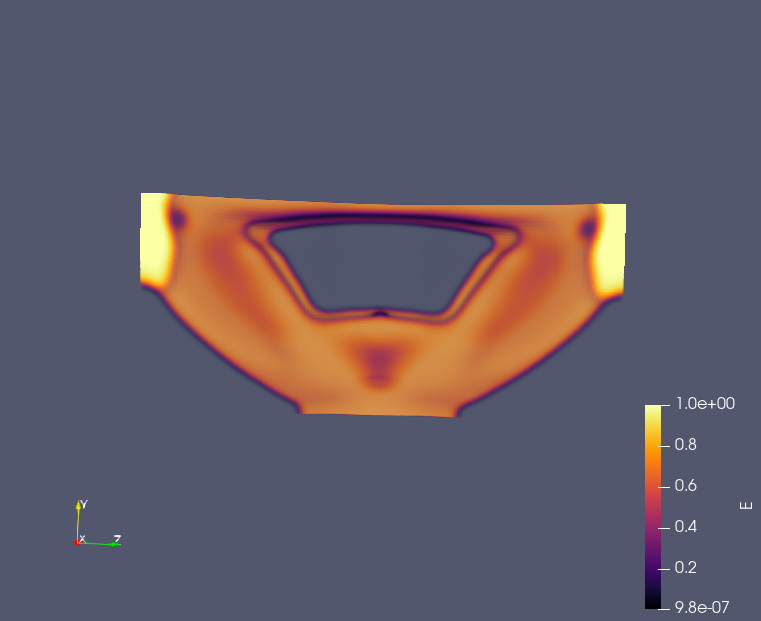}
        \caption{View 4.}
    \end{subfigure}
    \caption{Four views of the drone structure with 0.12 void volumetric fraction.}
    \label{fig:drone 012 fraction views}
\end{figure}

The computational domain is $1.0 \times 0.5 \times 1.0$ and the problem is solved on a $128 \times 64 \times 128$ structured grid.
For this problem, each optimization loop comprises 50 iterations of the accelerated pseudo-transition iterations, 50 iterations of the PT iterations, one update step of $\phi$s using sensitivity analysis, and one evolution step for $\phi$s using the Cahn-Hilliard model. 
The drone skeleton after \num{40000} optimization steps is presented in Figure~\ref{fig:drone 020 fraction views} for a void volumetric fraction of 0.2 and in Figure~\ref{fig:drone 012 fraction views} for a void volumetric fraction of 0.12.
The entire optimization process on a single GPU card completes in just 42 minutes, delivering in a highly refined structure for the drone skeleton that adheres closely to both the underlying physics and specific volumetric constraints. 

\FloatBarrier 

\FloatBarrier 

\section{Conclusion and Future Work}\label{conclusion}
In conclusion, we have introduced an innovative approach named PeTTO that combines the pseudo-transient and phase field methods to tackle topology optimization problems, leveraging the computational power of GPUs. 
By integrating PT and PF, we transformed the constrained optimization problem into solving a set of time-dependent PDEs, benefiting from insights from transient physics. 
The use of automatic differentiation streamlined sensitivity calculations, eliminating the error-prone manual derivations. Our hybrid approach, blending PT and accelerated PT methods, effectively balanced convergence speed and numerical stability, as demonstrate across various numerical examples.

Looking forward, we would also like to test this approach in topology optimization problems involving multi-physics which is common in diverse engineering disciplines. 
Furthermore, it would be interesting to explore more sophisticated constraints and integrate them into sensitivity calculations using automatic differentiation. 
Last, we aim to transition our approach to a multi-GPU environment and leverage JAX's multi-device parallel functionality to further accelerate topology optimization, advancing the frontier of high-performance computing in design optimization.

\newpage
\nocite{*}
\bibliography{main.bib}

\end{sloppypar}

\end{document}